\RequirePackage{fix-cm}
\documentclass[smallcondensed]{svjour3}     

\smartqed  
\usepackage{graphicx}
\usepackage{caption}
\usepackage{amsfonts}
\usepackage{mathptmx}

\usepackage[]{geometry}
\usepackage[square,numbers,sort&compress]{natbib}
\spnewtheorem{fact}{Fact}[section]{\bf}{\it}

\begin{document}
\title{On hypergraph cliques and polynomial programming
}


\author{Qingsong Tang         \and   Yuejian  Peng  \and   Xiangde Zhang \and Cheng Zhao
}


\institute{Qingsong Tang \at
              College of Sciences, Northeastern University, Shenyang, 110819, P.R.China. \at
             School of Mathematics, Jilin University, Changchun 130012, P.R. China.\\
              \email{t\_qsong@sina.com}
              \and
              Yuejian  Peng \at
              College of Mathematics, Hunan University, Changsha 410082, P.R. China. \at Supported  in part by National Natural Science Foundation of China (No. 11271116).\\
              \email{ypeng1@163.com}
               \and
           Xiangde Zhang (Corresponding author)\at
             College of Sciences, Northeastern University, Shenyang, 110819, P.R.China\\
             \email{zhangxdneu@163.com}
              \and
             Cheng Zhao \at
             Department of Mathematics and Computer Science, Indiana State University, Terre Haute, IN, 47809 USA. \at School of Mathematics, Jilin University, Changchun 130012, P.R. China.\\
             \email{cheng.zhao@indstate.edu}
}

\date{Received: date / Accepted: date}

\maketitle

\begin{abstract}
Motzkin and Straus established a close connection between the maximum
clique problem and a solution (namely graph-Lagrangians) to the maximum value of a class
of homogeneous quadratic multilinear functions over the standard simplex of the Euclidean
space in 1965. This connection provides a new proof of Tur\'an's theorem. Recently, an extension of Motzkin-Straus theorem was proved for non-uniform hypergraphs whose edges contain 1 or 2 vertices in \cite{PPTZ}. It is interesting if similar results hold for other non-uniform hypergraphs. In this paper, we give some connection between polynomial programming and the clique of non-uniform hypergraphs whose edges contain 1, or  2, and more vertices. Specifically, we obtain some Motzkin-Straus type results in terms of the graph-Lagrangian of non-uniform hypergraphs whose
edges contain 1, or 2, and more vertices.

\end{abstract}

\section{Introduction}
\label{intro}
In 1965,  Motzkin and Straus provided a new proof of Tur\'an's theorem based on a remarkable connection between the maximum clique and the graph-Lagrangian of a graph in \cite{MS}.  In fact, the connection of graph-Lagrangians and Tur\'{a}n densities can be used to give another proof of the fundamental theorem of Erd\H os-Stone-Simonovits on Tur\'{a}n densities of graphs. This type of connection aroused  interests in the study of graph-Lagrangians of $r$-graphs. A generalization of Motzkin-Straus theorem and Erd\H os-Stone-Simonovits theorem to non-uniform hypergraphs whose edges contain 1 or 2 vertices was given in \cite{PPTZ}.

A hypergraph $H=(V,E)$ consists of a vertex set $V$ and an edge set $E$, where every edge in $E$ is a subset of $V$. The set $T(H)=\{|F|:F \in E\}$ is called the set of  $edge\ types$ of $H$. We also say that $H$ is a $T(H)$-graph. For example, if $T(H)=\{1, 2\}$, then we say that $H$ is a $\{1, 2\}$-graph. If all edges have the same cardinality $r$, then $H$ is called an $r$-uniform hypergraph or $r$-graph. A $2$-uniform graph is called a graph. A hypergraph is non-uniform if it has at least two edge types.  For any $r \in T(H)$, the $level \ hypergraph \ H^r$ is the hypergraph consisting of all edges with $r$ vertices of $H$. We write $H^T_n$ for a hypergraph $H$ on $n$ vertices with $T(H)= T$. An edge $\{i_1, i_2, \cdots, i_r\}$ in a hypergraph  is simply written as  $i_1i_2\cdots i_r$ throughout the paper.

 For a positive integer $n$, let $[n]$ denote the set $\{1,2,\cdots,n\}$. For a finite set $V$ and a positive integer $i$, let ${V \choose i}$ denote the family of all $i$-subsets of $V$. The complete hypergraph $K^T_n$ is a hypergraph on $n$ vertices with edge set $ \bigcup_{i\in T} {[n] \choose i}$.   For example, $K_n^{\{r\}}$ is the complete $r$-uniform hypergraph on $n$ vertices.   $K_n^{[r]}$ is the non-uniform hypergraph with all possible edges of cardinality at most $r$. The complete graph on $n$ vertices $K_n^{\{2\}}$ is also called a clique. We also let $[t]^{(r)}$ represent the complete $r$-uniform hypergraph on vertex set $[t]$.

A useful tool in extremal problems of uniform hypergraphs (graphs) is the graph-Lagrangian of a uniform hypergraph (graph).
\begin{definition}
For  an $r$-uniform graph $H$ with the vertex set $[n]$, edge set $E(H)$, and a vector \\$\vec{x}=(x_1,\ldots,x_n) \in {\mathbb R}^n$,
we associate a homogeneous polynomial in $n$ variables, denoted by  $\lambda (G,\vec{x})$  as follows:
$$\lambda (H,\vec{x}):=\sum_{i_1i_2 \cdots i_r \in E(H)}x_{i_1}x_{i_2}\ldots x_{i_r}.$$
Let $S:=\{\vec{x}=(x_1,x_2,\ldots ,x_n): \sum_{i=1}^{n} x_i =1, x_i
\ge 0 {\rm \ for \ } i=1,2,\ldots , n \}$.
Let $\lambda (H)$ represent the maximum
 of the above homogeneous  multilinear polynomial of degree $r$ over the standard simplex $S$. Precisely
 $$\lambda (H): = \max \{\lambda (H, \vec{x}): \vec{x} \in S \}.$$
\end{definition}
The value $x_i$ is called the {\em weight} of the vertex $i$.
A vector $\vec{x}:=(x_1, x_2, \ldots, x_n) \in {\mathbb R}^n$ is called a feasible weighting for $G$ iff
$\vec{x}\in S$. A vector $\vec{y}\in S$ is called an {\em optimal weighting} for $G$
if $\lambda (G, \vec{y})=\lambda(G)$. We call $\lambda(G)$   the  graph-Lagrangian of $G$.



\begin{remark}
$\lambda (G)$ was called Lagrangian of $H$ in literature \cite{FF,PZ,FR84,T}. The terminology `graph-Lagrangian' was suggested by Franco Giannessi.
\end{remark}

Motzkin and Straus in \cite{MS} showed that  the graph-Lagrangian of a 2-graph is determined by the order of its maximum clique.

\begin{theorem} (\cite{MS}) \label{MStheo}
If $G$ is a 2-graph in which a largest clique has order $t$, then
$\lambda(G)=\lambda(K^{(2)}_t)=\lambda([t]^{(2)})={1 \over 2}(1 - {1 \over t})$.

Furthermore, the vector $\vec{x}=(x_1,x_2,\ldots ,x_n)$ given by $x_{i}={1 \over t}$ if $i$ is a  vertex in a fixed maximum complete $\{1,2\}$-subgraph and $x_i=0$ else is an optimal weighting.
\end{theorem}

This result provides a solution to the optimization problem for a class of homogeneous quadratic multilinear functions over the standard simplex of an Euclidean  plane.
The Motzkin-Straus result and its extension were successfully employed in optimization to provide heuristics for the
maximum clique problem \cite{B1,B2,B3,G9}.  It has  been also generalized to vertex-weighted graphs \cite{G9} and edge-weighted graphs with applications to pattern recognition in image analysis \cite{B1, B2, B3, G9, PP, PP15, RTP20}. An attempt to generalize the Motzkin-Straus theorem to hypergraphs is due to S\'os and
Straus\cite{SS}. Recently, in \cite{BP1,BP2} Rota Bul\'o and Pelillo generalized the Motzkin and Straus'
result to $r$-graphs in some way using a continuous characterization of maximal cliques
other than graph-Lagrangians of hypergraphs.

The graph-Lagrangian of a hypergraph has been a useful tool in hypergraph extremal problems.
For example, Sidorenko \cite{sidorenko89} and Frankl-Furedi \cite{FF} applied graph-Lagrangians of hypergraphs in
finding Tur\'an densities of hypergraphs. Frankl and R\"{o}dl \cite{FR84} applied it in disproving Erd\"os
long standing jumping constant conjecture. In most applications, we need an upper bound for the graph-Lagrangian of a
hypergraph.

Note that the graph-Lagrangian of an $r$-uniform graph can be viewed as the supremum of  densities of its blow-ups multiplying a constant (${1 \over r!}$). The graph-Lagrangian of a non-uniform hypergraph defined in \cite{PPTZ}  is the supremum of densities of its blow-ups.

\begin{definition}
For a hypergraph $H^{T}_n$ with $T(H)=T$ and a vector $\vec{x}=(x_1,\ldots,x_n) \in R^n$,
define
$$\lambda' (H^{T}_n,\vec{x}):=\sum_{r\in T}(r!\sum_{i_1i_2 \cdots i_r \in E(H^r)}x_{i_1}x_{i_2}\ldots x_{i_r}).$$

Let $S=\{\vec{x}=(x_1,x_2,\ldots ,x_n): \sum_{i=1}^{n} x_i =1, x_i \ge 0 {\rm \ for \ } i=1,2,\ldots , n \}$.
The Graph-Lagrangian of $H^{T}_n$, denoted by $\lambda' (H^{T}_n)$, is defined as
 $$\lambda' (H^{T}_n): = \max \{\lambda' (H^{T}_n, \vec{x}): \vec{x} \in S \}.$$
The value $x_i$ is called the {\em weight} of the vertex $i$. A vector $\vec{y}\in S$ is called an {\em optimal weighting} for $H$ if $\lambda' (H, \vec{y})=\lambda'(H)$.
\end{definition}
In \cite{PPTZ},  Peng et al. gave a generalization of Mozkin-Straus result to $\{1,2\}$-graphs.

\begin{theorem} (\cite{PPTZ}) \label{th1}
If $H$ is a $\{1,2\}$-graph and the order of its maximum complete $\{1,2\}$-subgraph is $t$, where $t\ge 2$, then
$\lambda '(H)=\lambda' (K^{\{1,2\}}_t)=2 - {1 \over t}$.

Furthermore, the vector $\vec{x}=(x_1,x_2,\ldots ,x_n)$ given by $x_{i}={1 \over t}$ if $i$ is a  vertex in a fixed maximum complete $\{1,2\}$-subgraph and $x_i=0$ else is an optimal weighting.
\end{theorem}



Some related Motzkin-Straus  type results in terms of graph-Lagrangians for non-uniform hypergraphs can be found in \cite{GLPS}.

In \cite{PP}, a more general question  is proposed.
\begin{problem}\label{generalproblem}
Let $H$ be  an $\{r_0,r_1,r_2,\ldots,r_m\}$-graph, $r_0<r_1<r_2< \ldots <r_m$, with vertex set $V(H)=[n]$ and
edge set $E(H)$.  Let $S=\{\vec{x}=(x_1,x_2,\ldots ,x_n)\in R^n: \sum_{i=1}^{n} x_i =1, x_i
\ge 0 {\rm \ for \ } i=1,2,\ldots , n \}$. Let $\alpha_{r_i}, 1\le i\le m$ be positive constants. For $\vec{x}\in S$,
let
\begin{eqnarray*}
L_{\{\alpha_{r_1},\alpha_{r_2},\cdots,\alpha_{r_m}\}} (H,\vec{x})&:=&\sum_{ i_1i_2 \ldots i_{r_0}  \in E(H^{r_0})}x_{i_1}x_{i_2}\ldots x_{i_{r_0}}+\alpha_{r_1} \sum_{ i_1i_2 \ldots i_{r_1}  \in E(H^{r_1})}x_{i_1}x_{i_2}\ldots x_{i_{r_1}}\\
&+&\ldots+\alpha_{r_m} \sum_{ i_1 i_2\ldots i_{r_m} \in E(H^{r_m})}x_{i_1}x_{i_2}\ldots x_{i_{r_m}}.
\end{eqnarray*}
The {\em polynomial optimization problem} of $H$ is
\begin{eqnarray}\label{lg}
L_{\{\alpha_{r_1},\alpha_{r_2},\cdots,\alpha_{r_m}\}}  (H): = \max \{L(H, \vec{x}): \vec{x} \in S \}.
\end{eqnarray}
We sometimes simply write $L_{\{\alpha_{r_1},\alpha_{r_2},\cdots,\alpha_{r_m}\}} (H,\vec{x})$ and $L_{\{\alpha_{r_1},\alpha_{r_2},\cdots,\alpha_{r_m}\}}(H)$ as $L(H,\vec x)$ and $L(H)$ if there is no confusion.
The value $x_i$ is called the {\em weight} of the vertex $i$. A vector $\vec{x}=(x_1, x_2, \ldots, x_n) \in {\mathbb R}^n$ is called a feasible solution to (\ref{lg}) if and only if $\vec{x}\in S$. A vector $\vec{y}\in S$ is called a {\em solution} to optimization problem (\ref{lg}) if and only if $L (H, \vec{y})=L(H)$.
\end{problem}

\begin{remark}\label{weightedlagrangians}
Let $H$ be  an $\{r_0,r_1,r_2,\ldots,r_m\}$-graph, $r_0<r_1<r_2< \ldots <r_m$, with vertex set $V(H)=[n]$ and
edge set $E(H)$. Clearly, $\lambda' (H,\vec{x})=r_0!L_{\{\frac{r_1!}{r_0!},\frac{r_2!}{r_0!},\cdots,\frac{r_m!}{r_0!}\}}  (H, \vec{x}).$ Hence we can view $L(H)$ as subgraph weighted graph-Lagrangian of $H$.
\end{remark}

Peng etc. in \cite{PP} gave some Motzkin-Straus type results to $\{1,r \}$-graphs and  $\{1, 2, 3\}$-graphs  for the polynomial programming (\ref{generalproblem}).

\begin{theorem} \cite{PP} \label{thm4}
Let $\alpha_r>0$ be a constant. Let $H$ be a $\{1,r\}$-graph. If both the order of its maximum complete $\{1,r\}$-subgraph  and  the order of its maximum complete $\{1\}$-subgraph are $t$, where $\displaystyle {t\geq \lceil {[\alpha_r-(r-2)!]^{r-2} \over (r-2)!\alpha_r^{r-3}}\rceil }$, then \[L _{\{ \alpha_r \}}(H) = L_{\{ \alpha_r \}} \left( {{K_t}^{\{ 1,r\} }} \right) ={1+ \alpha_r \frac{\prod_{i=1}^{r-1} (t-i)}{r!t^{r-1}}}.\]
Furthermore, the vector $\vec{x}=(x_1,x_2,\ldots ,x_n)$ given by $x_{i}={1 \over t}$ if $i$ is a  vertex in a fixed maximum complete $\{1,r\}$-subgraph and $x_i=0$ else is a solution to the optimization problem (\ref{lg}) with $m=1$ and $r_0=1$.
\end{theorem}

\begin{theorem} \cite{PP} \label{th4}
Let $\alpha_2,\alpha_3>0$ be constants. Let $H$ be a $\{1, 2, 3\}$-graph. If both the order of its maximum complete
$\{1,2, 3\}$-subgraph  and  the order of its maximum complete
$\{1\}$-subgraph are $t$, where $t\ge \lceil{ (\alpha_2+\alpha_3)^2-\alpha_3 \over \alpha_2+\alpha_3}\rceil $, then
\[ L_{\{ \alpha_2,\alpha_3 \}}(H) =  L_{\{ \alpha_2, \alpha_3 \}} \left( {{K_t}^{\{ 1, 2, 3\} }} \right) ={1+  \alpha_2 \frac{t-1}{2t}+  \alpha_3 \frac{(t-1)(t-2)}{6t^2}}.\]
Furthermore, the vector $\vec{x}=(x_1,x_2,\ldots ,x_n)$ given by $x_{i}={1 \over t}$ if $i$ is a  vertex in a fixed maximum complete $\{1,2,3\}$-subgraph and $x_i=0$ else is a solution to the corresponding optimization problem.
\end{theorem}

In this paper, we will prove other Motzkin-Straus  type results to non-uniform hypergraphs whose edges contain 1, $2$, and more vertices for (\ref{lg}). Here are our main results.
\begin{theorem} \label{thL1}
(a) Let $\alpha_r>0$ be a constant. Let $H$ be a $\{2,r\}$-graph. If both the order of its maximum complete $\{2,r\}$-subgraphs and the vertex order of $H^2$ are $t$, where $t\geq  \frac{\alpha_r}{(r-2)!}+1$, then
\[L _{\{ \alpha_r \}}(H) = L_{\{ \alpha_r \}} \left( {{K_t}^{\{ 2,r\} }} \right) ={\frac{t-1}{2t}+ \alpha_r \frac{\prod_{i=1}^{r-1} (t-i)}{r!t^{r-1}}}.\]

 (b) Let $\alpha_2,\alpha_r>0$ be constants.  Let $H$ be a $\{1,2,r\}$-graph. If both the order of its maximum complete $\{1,2,r\}$-subgraphs and the vertex order of $H^2$ are $t$, where $t\geq  \frac{\alpha_r}{\alpha_2(r-2)!}+1$, then
\[ L_{\{ \alpha_2,\alpha_r \}}(H) =  L_{\{ \alpha_2, \alpha_r \}} \left( {{K_t}^{\{ 1, 2, r\} }} \right) ={1+  \alpha_2 \frac{t-1}{2t}+\alpha_r \frac{\prod_{i=1}^{r-1} (t-i)}{r!t^{r-1}}}.\]

Furthermore, the vector $\vec{x}=(x_1,x_2,\ldots ,x_n)$ given by $x_{i}={1 \over t}$ if $i$ is a  vertex in a fixed maximum complete $\{1,2,3\}$-subgraph and $x_i=0$ else is a solution to the corresponding optimization problem in  both (a) and (b).
\end{theorem}

\begin{theorem} \label{thL2}
(a) Let $\alpha_r>0$ be a constant. Let $H$ be a $\{2,r\}$-graph. If  the order of its maximum complete $\{2,r\}$-subgraphs is $t$,  and the number of edges in $H^2$, say $m$, satisfies ${t \choose 2} \le m \le {t \choose 2}+t-2$, where $t\geq  \frac{\alpha_r}{(r-2)!}+1$, 
then
\[L _{\{ \alpha_r \}}(H) = L_{\{ \alpha_r \}} \left( {{K_t}^{\{ 2,r\} }} \right) ={\frac{t-1}{2t}+ \alpha_r \frac{\prod_{i=1}^{r-1} (t-i)}{r!t^{r-1}}}.\]

 (b) Let $\alpha_2,\alpha_r>0$ be constants satisfying $\alpha_2\geq\frac{\alpha_r}{(r-2)!}$.  Let $H$ be a $\{1,2,r\}$-graph. If the order of its maximum complete $\{1,2,r\}$-subgraphs is $t$,  and the number of edges in $H^2$, say $m$, satisfies ${t \choose 2} \le m \le {t \choose 2}+t-2$, where $t\geq  \frac{\alpha_r}{\alpha_2(r-2)!}+1$ and $  \alpha_2\geq\frac{ \alpha_r}{2(r-2)!}$, then
\[ L_{\{ \alpha_2,\alpha_r \}}(H) =  L_{\{ \alpha_2, \alpha_r \}} \left( {{K_t}^{\{ 1, 2, r\} }} \right) ={1+  \alpha_2 \frac{t-1}{2t}+  \alpha_r \frac{\prod_{i=1}^{r-1} (t-i)}{r!t^{r-1}}}.\]

Furthermore, the vector $\vec{x}=(x_1,x_2,\ldots ,x_n)$ given by $x_{i}={1 \over t}$ if $i$ is a  vertex in a fixed maximum complete $\{1,2,3\}$-subgraph and $x_i=0$ else is a solution to the corresponding optimization problem in  both (a) and (b).
\end{theorem}
Applying Theorems \ref{thL1}, \ref{thL2}, Remark \ref{weightedlagrangians}, and by choosing appropriate coefficients in the polynomial programming (\ref{generalproblem}), it is easy to see that the following results hold.
\begin{corollary} \label{th6}
(a) Let $H$ be a $\{2,r\}$-graph. If both the order of its maximum complete $\{2,r\}$-subgraphs and the vertex order of $H^2$ are $t$, where $t\geq  \frac{r(r-1)}{2}+1$, then
$\lambda '(H)=\lambda '(K^{\{2,r\}}_{t})$.

 (b) Let $H$ be a $\{1,2,r\}$-graph. If both the order of its maximum complete $\{1,2,r\}$-subgraphs and the vertex order of $H^2$ are $t$, where $t\geq  \frac{r(r-1)}{2}+1$, then
$\lambda '(H)=\lambda '(K^{\{1,2,r\}}_{t})$.

Furthermore, the vector $\vec{x}=(x_1,x_2,\ldots ,x_n)$ given by $x_{i}={1 \over t}$ if $i$ is a  vertex in a fixed maximum complete $\{1,2,3\}$-subgraph and $x_i=0$ else is a solution to the corresponding optimization problem in  both (a) and (b).
\end{corollary}

\begin{corollary} \label{thlast}
(a)  Let $3\leq r\leq 4.$ Let $H$ be a $\{2,r\}$-graph. If  the order of its maximum complete $\{2,r\}$-subgraphs is $t$,  and the number of edges in $H^2$, say $m$, satisfies ${t \choose 2} \le m \le {t \choose 2}+t-2$, where $t\geq \frac{r(r-1)}{2}+1$, then
$$\lambda '(H)=\lambda '(K^{\{2,r\}}_{t})={\frac{t-1}{t}+   \frac{\prod_{i=1}^{r-1} (t-i)}{t^{r-1}}}.$$
(b)  Let $3\leq r\leq 4.$ Let $H$ be a $\{1,2,r\}$-graph. If the order of its maximum complete $\{1,2,r\}$-subgraphs is $t$,  and the number of edges in $H^2$, say $m$, satisfies ${t \choose 2} \le m \le {t \choose 2}+t-2$, where $t\geq \frac{r(r-1)}{2}+1$, then
$$\lambda '(H)=\lambda '(K^{\{1,2,r\}}_{t})={1+   \frac{t-1}{t}+   \frac{\prod_{i=1}^{r-1} (t-i)}{t^{r-1}}}.$$

Furthermore, the vector $\vec{x}=(x_1,x_2,\ldots ,x_n)$ given by $x_{i}={1 \over t}$ if $i$ is a  vertex in a fixed maximum complete $\{1,2,3\}$-subgraph and $x_i=0$ else is a solution to the corresponding optimization problem in both (a) and (b).
\end{corollary}

The rest of the paper is organized as follows. Some useful results are summarized in Section \ref{sec:1}. The proofs of
Theorems \ref{thL1}, \ref{thL2} are given in Section 3. Further Motzkin-Straus type results for $\{2, r_3, ..., r_m\}$-graphs and $\{1,2, r_3, ..., r_m\}$-graphs are given in Section 3 as well.

\section{Some Preliminary  Results}
\label{sec:1}
 We will impose an additional condition on any solution $\vec{x}=(x_1,x_2,\ldots ,x_n)$ to the polynomial programming (\ref{generalproblem}):

(i)$x_1 \ge x_2 \ge \ldots \ge x_n \geq 0$.

(ii)$ \vert\{i : x_i > 0 \}\vert$ is minimal, i.e., if $\vec y $ is a feasible solution to the polynomial programming (\ref{generalproblem}) satisfying $|\{i : y_i > 0\}| < |\{i : x_i > 0\}|$, then $ L (H, {\vec y}) < L(H)$ .

For a hypergraph $H=(V,E)$, $i\in V$, and $r\in T(H)$, let $E_i^r=\{A \in V^{(r-1)},  A \cup \{i\} \in E^r\}$. For a pair of vertices $i,j \in V$, let $E_{ij}^r=\{B \in V^{(r-2)} B \cup \{i,j\} \in E^r\}$.
Let
$(E^{r}_i)^c=\{A \in V^{(r-1)}, A \cup \{i\} \in V^{(r)} \backslash E\}$, $(E^r_{ij})^c=\{B \in V^{(r-2)} B \cup \{i,j\} \in V^{(r)} \backslash E^r\}$, and $E_{i\setminus j}^r=E_i^r\cap (E^{r}_{j})^c.$ Let $L(E_i^r, \vec{x})=\alpha_{r}\lambda (E_i^r, \vec{x})$, where $\alpha_{r_0}=1$.  And $L(E_{ij}^r, \vec{x})$ and $L(E_{i\setminus j}^r, \vec{x})$ are defined similarly.

Let $E_i=\cup_{r\in T(H)} E_i^r$, $E_{i\setminus j}=\cup_{r\in T(H)} E_{i\setminus j}^r$, and   $E_{ij}=\cup_{r\in T(H)} E_{ij}^r$. Let  $L(E_i, \vec{x})=\cup_{r\in T(H)} L(E_i^r, \vec{x})$.  And $L(E_{ij}, \vec{x})$ and $L(E_{i\setminus j}, \vec{x})$ are defined similarly. Note that $L(E_i, \vec{x})={\partial L (H, {\vec x})\over \partial x_i}$ and $L(E_{ij}, \vec{x})={\partial L (H, {\vec x})\over \partial x_i\partial x_j}$.

Let $H=([n],E)$. For $e\in E$, and $i,j\in [n]$ with $i<j$, define
\begin{eqnarray*}
C_{i\leftarrow j}(e) :&=& \left\{
\begin{array}{rl}
(e\backslash\{j\})\cup \{i\}&   \  {\rm if}  \ i\notin e\ {\rm and}  \ j\in e,\\
e& \ {\rm otherwise}.
\end{array}\right.
\end{eqnarray*}

and $\mathcal{C}_{i\leftarrow j}(e)=\{C_{i\leftarrow j}(e):e\in E\}\bigcup \{e,C_{i\leftarrow j}(e)\in E\}.$

We say that $H$ is left-compressed if $C_{i\leftarrow j}(E)=E$ for every $1\leq i\leq j.$
\begin{remark}(Equivalent definition of  left-compressed)
 A  $T(H)$-hypergraph $H=([n],E)$ is \emph{left-compressed} if and only if for any $r \in T(H)$, $j_1j_2 \cdots j_r \in E$ implies $i_1i_2\cdots i_r \in E$ provided $i_p \le j_p$ for every $p$, $1 \le p\le r$. Equivalently, a $T(H)$-hypergraph $H=([n],E)$ is \emph{left-compressed} if and only if for any $r \in T(H)$, $E^r_{j\backslash i}=\emptyset$ for any $1 \le i< j\le n.$
\end{remark}

\begin{lemma} (\cite{PP}) \label{leftcom}
Let $H=([n],E)$ be a $T(H)$-graph,
$i,j\in [n]$ with $i<j$ and $\vec x = (x_1,\cdots, x_n)$ be a solution to the polynomial programming (\ref{generalproblem}). Write $H_{i\leftarrow j}=([n],\mathcal{C}_{i\leftarrow j}(E))$. Then,
$$L(H, \vec x)\leq L(H_{i\leftarrow j}, \vec x).$$
\end{lemma}

\begin{lemma}\cite{PP} \label{Lemmkkt}
If $x_1 \ge x_2 \ge \ldots \ge x_k >x_{k+1}=x_{k+2}=\ldots =x_n=0$ and ${\vec x}=(x_1, x_2, \ldots, x_n)$ is a solution to  the polynomial programming (\ref{generalproblem}), then (a) ${\partial L (H, {\vec x})\over \partial x_1}={\partial  L (H, {\vec x})\over \partial x_2}= \ldots ={\partial  L (H, {\vec x})\over \partial x_k}.$  This is equivalent to  $L(E_i, {\vec x})=L(E_j, {\vec x})$ for $1\le i<j\le k$.
(b) $\forall 1\le i<j\le k,$ there exists an edge $e \in E(H)$  such that $\{i,j\} \subseteq e$.
\end{lemma}

\begin{remark}\label{r1} (a) Lemma \ref{Lemmkkt} part (a) implies that
$$x_jL(E_{ij}, {\vec x})+L (E_{i\setminus j}, {\vec x})=x_iL(E_{ij}, {\vec x})+L(E_{j\setminus i}, {\vec x}).$$
In particular, if $H$ is left-compressed, then
$$(x_i-x_j)L(E_{ij}, {\vec x})=L(E_{i\setminus j}, {\vec x})$$
for any $i, j$ satisfying $1\le i<j\le k$ since $E_{j\setminus i}=\emptyset$.

(b) If  $G$ is left-compressed, then for any $i, j$ satisfying $1\le i<j\le k$,
\begin{equation}\label{enbhd}
x_i-x_j={L (E_{i\setminus j}, {\vec x}) \over L(E_{ij}, {\vec x})}
\end{equation}
holds.  If  $G$ is left-compressed and  $E_{i\setminus j}=\emptyset$ for $i, j$ satisfying $1\le i<j\le k$, then $x_i=x_j$.

(c) By (\ref{enbhd}), if  $H$ is left-compressed, then a solution   ${\vec x}=(x_1, x_2, \ldots, x_n)$ to the optimization problem (\ref{lg})  must satisfy
\begin{equation}\label{conditiona}
x_1 \ge x_2 \ge \ldots \ge x_n \ge 0.
\end{equation}
\end{remark}

In \cite{PZ}, \cite{TPZZ2}, and \cite{PTZ}, the following theorems for $3$-graphs and $r$-graphs were proved, respectively.

\begin{theorem}  \cite{PZ} \label{theoPZ}
 Let $m$ and $t$ be positive integers satisfying $${t \choose 3} \le m \le {t \choose 3} + {t-1 \choose 2}.$$  Let $H$ be a $3$-graph  with $m$ edges and containing a clique of order $t$. Then $\lambda(H) = \lambda([t]^{(3)})$.
\end{theorem}
\begin{theorem}  \cite{TPZZ2} \label{theoTPZZ} Let $m$ and $t$ be integers satisfying
${t \choose 3} \le m \le {t \choose 3} + {t-1 \choose 2} - \frac{t}{2}.$
Let $G$ be a $3$-graph with $m$ edges, if $G$ does not contain a complete subgraph of order  $t$, then  $\lambda(G)<\lambda([t]^{(3)}).$
\end{theorem}
\begin{theorem}  \cite{PTZ} \label{theoPTZ}
 Let $m$ and $t$ be positive integers satisfying $${t \choose r} \le m \le {t \choose r} + {t-1 \choose r-1}-(2^{r-3}-1)({t-1 \choose r-2}-1).$$  Let $H$ be an $r$-graph on $t+1$ vertices with $m$ edges and containing a clique of order $t$. Then $\lambda(G) = \lambda([t]^{(r)})$.
\end{theorem}

\section{Proofs of main results}
\label{sec:3}

In order to prove Theorems \ref{thL1} and \ref{thL2}, we begin with two lemmas. In the rest of the paper£¬  an optimal (feasible)  weighting for $H$ refers to a  solution (feasible) to the polynomial programming (\ref{generalproblem}) unless specifically stated.
\begin{lemma} \label{th7}
 (a) Let $H$ be a $\{2,r\}$-graph. If both the order of its maximum complete $\{2,r\}$-subgraphs and the vertex order of $H^2$ are $t$, and $H^r$ is $[s]^{(r)}$,  where $s\geq t\geq \frac{\alpha_r}{\alpha_2(r-2)!}+1$,  then
$L(H)=L(K^{\{2,r\}}_{t})$.

(b) Let $H$ be a $\{1,2,r\}$-graph. If both the order of its maximum complete $\{1,2,r\}$-subgraphs and  the vertex order of $H^2$ are $t$, $H^1$ is $[u]^{(1)}$, and $H^r$ is $[v]^{(r)}$, where $u\geq 4$ and $v\geq 4$, then
$L(H)=L(K^{\{1,2,r\}}_{t})$.


\end {lemma}


\noindent{\em Proof of Lemma \ref{th7}}  We only give the proof of (b). The proof of (a) is similar to (b).

Applying Lemma \ref{Lemmkkt}(a) and a direct calculation, we get a solution $\vec{y}$ to for the polynomial programming (\ref{generalproblem}) when $H={K_t}^{\{ 1,r\} }$ which is given by $y_i=1/t$ for each $i(1\le i\le t)$ and $y_i=0$ else. So $ L({{K_t}^{\{ 1,2,r\} }}) ={1+ \alpha_2\frac{t-1}{2t}+\alpha_r \frac{\prod_{i=1}^{r-1} (t-i)}{r!t^{r-1}}}$. Since $K^{\{1,2,r\}}_{t}\subset H$, clearly $L(H)\geq L(K^{\{1,2,r\}}_{t}).$ Thus to prove Thorem \ref{th7}, we only need to prove that  $L(H)\leq L(K^{\{1,2,r\}}_{t})$.

Denote $M_{(s,t,\{1,2,r\})}=\max\{L(H): H $ is a $\{1,2,r\}$-graph with $H^1=[u]^{(1)}$, $H^r=[v]^{(r)}$ and both the order of its maximum complete $\{1, 2,r\}$-subgraph and the (vertex) order of $H^2$ are $t$\}. We can assume that $L(H)=M(t+1,t,\{1,2,r\})$,i.e. $H$ is an extremal graph.  We can assume that $H$ is left-compressed. If $H$ is not left-compressed, performing a sequence of left-compressing operations(i.e. replace $E$ by
$\mathcal{C}_{i\leftarrow j}(E)$ if $\mathcal{C}_{i\leftarrow j}(E)\neq E$), we will get  a left-compressed $\{1,2,r\}$-graph with $H'$ with the same number of edges,  $H'^1=[u]^{(1)}$, $H'^r=[v]^{(r)}$, and both the order of its maximum complete $\{1,2,r\}$-subgraph and the (vertex) order of $H'^{2}$ are still $t$. By Lemma \ref{leftcom} $H'$ is also an extremal graph. We give the proof of the case  $u\leq v$ below. The proof for the case $v\le u$ is similar and the details will not be given.

Let $\vec{x}=(x_{1},x_{2},\ldots ,x_{v})$ be an optimal weighting of $H$. By Remark \ref{r1}(c),  $x_1 \ge x_2 \ge \ldots  \geq x_{v}\geq 0$.
By Remark \ref{r1}(b) we may assume  that $x_1=\cdots =x_t$, $x_{t+1}=\cdots=x_{u}$ and $x_{u+1}=\cdots=x_{v}$.

First we show that $x_v=0$. Assume that $x_v>0$ for a contradiction.  By Lemma \ref{Lemmkkt}(a), $L (E_1, {\vec x})=L (E_v, {\vec x})$. Assume $u<v$ since otherwise we only need to prove that  $x_u=0$. Hence
\begin{eqnarray}\label{eqcit1}
1+\alpha_2(1-x_1-x_{t+1}-\cdots -x_v)+L(E_{1\backslash v}^{r},\vec{x})+x_{v}L(E_{1v}^{r},\vec{x})-x_{1}L(E_{1v}^{r},\vec{x})=0.
\end{eqnarray}
Since $1-x_1-x_{t+1}-\cdots -x_v=(t-1)x_1$, 
the above equality is equivalent to
\begin{eqnarray*}
\frac{1}{\alpha_2}+(t-1)x_1=\frac{x_1-x_v}{\alpha_2}L(E_{1v}^{r},\vec{x}).
\end{eqnarray*}
Since $0<L(E_{1v}^{r},\vec{x})\leq \frac{\alpha_r(1-x_1-x_{v})^{r-2}}{(r-2)!}<\frac{\alpha_r}{(r-2)!}$, then $(t-1)x_1<\frac{\alpha_r}{\alpha_2(r-2)!}x_1,$ i.e. $t<\frac{\alpha_r}{\alpha_2(r-2)!}+1$, which contradicts to $t\geq \frac{\alpha_r}{\alpha_2(r-2)!}+1$.

Since $x_{u+1}=\cdots=x_{v}=0$, we can assume that the hypergraph is on $[u]$.  Next we show that $x_u=0$. Assume that $x_u>0$ for a contradiction.  By Lemma \ref{Lemmkkt}, $L (E_1, {\vec x})=L (E_u, {\vec x})$. Hence
\begin{eqnarray}\label{eqcit2}
\alpha_2(1-x_1-x_{t+1}-\cdots -x_u)+L(E_{1\backslash u}^{r},\vec{x})+x_{u}L(E_{1u}^{r},\vec{x})-x_{1}L(E_{1u}^{r},\vec{x})=0.
\end{eqnarray}
Since $0<L(E_{1u}^{r},\vec{x})\leq \alpha_r\frac{(1-x_1-x_{u})^{r-2}}{(r-2)!}<\frac{\alpha_r}{(r-2)!}.$ Hence $(t-1)x_1<\frac{\alpha_r}{\alpha_2(r-2)!}x_1,$ i.e. $t<\frac{\alpha_r}{\alpha_2(r-2)!}+1$, which contradicts to $t\geq\frac{\alpha_r}{\alpha_2(r-2)!}+1$. This completes the proof of (b).\qed

\noindent{\em Proof of Theorem \ref{thL1} }  We only give the proof of (b). The proof of (a) is similar to (b). Since $K^{\{1,2,r\}}_{t}\subset H$, clearly $L(H)\geq L(K^{\{1,2,r\}}_{t}).$ Thus we only need to prove that  $L(H)\leq L(K^{\{1,2,r\}}_{t})$. Assume $H^{r}$ is on vertex set $[n]$. Let $\overline{E}^{1}=[n]^{(1)}$, $\overline{E}^{r}=[n]^{(r)}$ and $\overline{H}=E^{2}\bigcup \overline{E}^{1}\bigcup \overline{E}^{r}$. Then $L(H)\leq L(\overline{H})$. By Lemma \ref{th7} (b), $L(H)\le L(K^{\{1,2,r\}}_{t})$. This completes the proof of (b).\qed

Now we are ready to prove Theorem \ref{thL2} by applying Theorem \ref{thL1}.

\noindent{\em Proof of Theorem \ref{thL2}} We only give the proof of (b). The proof of (a) is similar to (b).
 Since $K^{\{1,2,r\}}_{t}\subset H$, clearly $L(H)\geq L(K^{\{1,2,r\}}_{t}).$ Thus to prove Theorem \ref{thlast}, we only need to prove that  $L(H)\leq L(K^{\{1,2,r\}}_{t})$.

Denote $M_{(m,t,\{1,2,r\})}=\max\{L(H): H $ is a $\{1,2,r\}$-graph with $m$ edges in $H^2$.   The order of its maximum complete $\{1,2,r\}$-subgraphs is $t\}$. We can assume that $L(H)=L{(m, t,\{1,2,r\})}$, i.e., $H$ is an extremal graph.  We can assume that $H$ is left-compressed. If $H$ is not left-compressed, performing a sequence of left-compressing operations (i.e. replacing $E$ by $\mathcal{C}_{i\leftarrow j}(E)$ if $\mathcal{C}_{i\leftarrow j}(E)\neq E$), we will get  a left-compressed $\{1,2,r\}$-graph $H'$ with the same number of edges.  And the order of its maximum complete $\{1,2,r\}$-subgraphs is  still $t$. By Lemma \ref{leftcom} $H'$ is also an extremal graph.

Let $\vec{x}=(x_{1},x_{2},\ldots ,x_{n})$ be an optimal weighting for $H$. Then $x_1 \ge x_2\geq\ldots \ge x_k>x_{k+1}= \ldots = x_{n}= 0$. First we show that the vertex order of $H^2$ is at most $t+1$.  Assume that there is  $ij\in H^2$ with $k\geq j\geq t+2.$   We define a new feasible weighting ${\vec y}$ for $H$ as follows. Let $y_l=x_l$ for $l\neq j-1, j$, $y_{j}=0 $ and $y_{j-1}=x_{j-1}+x_{j}$.  By Lemma \ref{Lemmkkt}(a), we have $ L (E_{j-1}, {\vec x})=L (E_j, {\vec x})$. Note that $(j-1)j\notin E^2$ for $j\geq t+1$. Hence
\begin{eqnarray*}\label{eqthlast}
L(H,\vec {y})- L(H,\vec {x})&=&\sum\limits_{u\in \{1, 2, r\}}x_{j}[L(E_{j-1}^u, \vec{x})-L(E_{j}^u, \vec{x})]-\sum\limits_{u\in\{1,2,r\}} x_{j}^2 L(E_{(j-1)j}^{u},\vec{x}) \nonumber \\
&=&-x_{j}^2 L(E_{(j-1)j}^{r},\vec{x}).
\end{eqnarray*}
Since $y_{j}=0$, we may remove all the edges containing $j$ from $E$ to form a new $3$-graph $\overline{H}=([n], \overline{E})$ with$\vert \overline{E}\vert=\vert E\vert-\vert E_{j}\vert$ and $L(\overline{H},\vec {y})=
L(H,\vec {y})$. Since $m \le {t \choose 2}+t-2$, we have  $(t-1)(j-1)\notin E^2$. Let $\overline{\overline{H}}=\overline{H}\bigcup \{(t-1)(j-1)\}$. Then $\overline{\overline{H}}$ is a $\{1,2,r\}$-graph. The order of its maximum complete $\{1,2,r\}$-subgraph is still $t$.  The number of edges in $\overline{\overline{H}}^2$ satisfies ${t \choose 2} \le m \le {t \choose 2}+t-2$.  Recalling that $\alpha_2\geq\frac{ \alpha_r}{2(r-2)!}$ and $x_1 \ge x_2\geq\ldots \ge x_{t-1}\geq x_{j-1}\geq x_{j}>0,$ we have
\begin{eqnarray*}\label{eqthlast}
L(\overline{\overline{H}},\vec {y})-L(H,\vec {x})&=& \alpha_2x_{t-1}(x_{j-1}+x_{j})-x_{j}^2L(E_{(j-1)j}^{r},\vec{x})\\
&>&\alpha_2x_{t-1}(x_{j-1}+x_{j})-\frac{\alpha_r}{(r-2)!}x_{j}^2\\
&\geq&0
\end{eqnarray*}
since $0<L(E_{(j-1)j}^{r},\vec{x})\leq \frac{\alpha_r(1-x_{j-1}-x_{j})^{r-2}}{(r-2)!}<\frac{\alpha_r}{(r-2)!}.$ This contradicts to that $H$ is an extremal graph. Hence the order of the $2$-graph is at most $t+1$.

Next we prove $L(H)\leq L(K^{\{1,2,r\}}_{t})$.

Let $\vec{x}=(x_{1},x_{2},\ldots ,x_{n})$ be an optimal weighting for $H$. Then $x_1 \ge x_2\geq\ldots \ge x_k>x_{k+1}= \ldots = x_{n}= 0$. We define a new feasible weighting ${\vec z}$ as follows. Let $z_l=x_l$ for $l\neq t, t+1$, $z_{t}=0 $ and $z_{t+1}=x_{t}+x_{t+1}$.  By Lemma \ref{Lemmkkt}(a), we have $L (E_t, {\vec x})=L (E_{t+1}, {\vec x})$. Note that $t(t+1)\notin E^2.$ Hence
\begin{eqnarray*}\label{eqthlast}
L(H,\vec {z})- L(H,\vec {x})&=&\sum\limits_{u\in \{1,2,r\}}x_{t}[L(E_{k}^u, \vec{x})-L(E_{t}^u, \vec{x})]-\sum\limits_{u\in\{2,r\}}x_{t}^2L(E_{t(t+1)}^{u},\vec{x}) \nonumber \\
&=&-x_{t}^2L(E_{t(t+1)}^{r},\vec{x}).
\end{eqnarray*}

Since $z_{t}=0$ we may remove all the edges containing $t$ from $E$ to form a new $3$-graph $H^*=([k], E^*)$ with$\vert E^*\vert=\vert E\vert-\vert E_{t}\vert$ and $L(H^*,\vec {y})=
L(H,\vec {y})$. Since $m \le {t \choose 2}+t-2$, we have  $(t-1)(t+1)\notin E^2$. Let $H^{**}:=H^*\bigcup \{(t-1)(t+1)\}$. Then $H^{**}$ is a $\{1,2,r\}$-graph.  Recalling that $\alpha_2\geq\frac{ \alpha_r}{2(r-2)!}$ and $x_1 \ge x_2\geq\ldots \ge x_{t-1}\geq x_{t}\geq 0,$   we have
\begin{eqnarray*}\label{eqthlast}
L(H^{**},\vec {y})- L(H,\vec {x})&=& \alpha_2x_{t-1}(x_{t}+x_{t+1})-x_{t}^2\lambda(E_{t(t+1)}^{r},\vec{x})\\
&\geq &\alpha_2x_{t-1}(x_{t}+x_{t+1})-\frac{\alpha_r}{(r-2)!}x_{t}^2\\
&\geq &0.
\end{eqnarray*}
since $0<L(E_{t(t+1)}^{r},\vec{x})\leq \frac{\alpha_r(1-x_{t}-x_{t+1})^{r-2}}{(r-2)!}<\frac{\alpha_r}{(r-2)!}.$ Since the  vertex order of  $(H^{**})^{2}$ is $t$, we have $L(H^{**})\le L(K_{t}^{\{1,2,r)\}})$ by Theorem \ref{thL1}. Hence $L(H)=L(H,\vec{x})\leq L(H^{**},\vec{y})\leq L(H^{**})=L(K_{t}^{\{1,2,r\}})$. This completes the proof of part (b).
 \qed

Using the method given in the proof of Theorem \ref{thL1}, we may generalize Theorem \ref{thL1} for $\{2,r_3,\cdots, r_l\}$-graphs ($\{1,2,r_3,\cdots, r_l\}$-graphs, respectively) in the following way,  where $r_3< \ldots <r_m$, with vertex set $V(H)=[n]$ and
edge set $E(H)$.

\begin{theorem}\label{Theogeneral}
 (a) Let $H$ be a $\{2,r_3,\cdots, r_m\}$-graph,  where $r_3< \ldots <r_m$. If both the order of its maximum complete $\{2,r_3,\cdots, r_m\}$-subgraph and the order of $\{2\}$-graph are $t$,  where $t\geq (m-2)\frac{\alpha_{r_l}}{\alpha_2(r_l-2)!}+1$ then
$$L(H)=L(K^{\{2,r_3,\cdots, r_m\}}_{t}).$$

(b) Let $H$ be a $\{1,2,r_3,\cdots, r_m\}$-graph,  where $r_3< \ldots <r_m$. If both the order of its maximum complete $\{1,2,r_3,\cdots, r_m\}$-subgraph and the order of $\{2\}$-graph are $t$,  where $t\geq (m-2) \frac{\alpha_{r_l}}{\alpha_2(r-2)!}+1$ then
$$L(H)=L(K^{\{1,2,r_3,\cdots, r_m\}}_{t}).$$


\end {theorem}
We remark that the proof of Theorem \ref{Theogeneral} is similar to the proof of Theorem \ref{thL1}. For instance, to prove Theorem \ref{Theogeneral}(b),  we  change (\ref{eqcit1}) to
\begin{eqnarray*}
1&+&\alpha_2(1-x_1-x_{t+1}-\cdots -x_v)+\sum\limits_{i\in \{3,\ldots, m\}}L(E_{1\backslash v}^{r_i},\vec{x})+\sum\limits_{v\in \{3,\ldots ,m\}}x_{v}L(E_{1v}^{r_i},\vec{x})\\
&-&\sum\limits_{i\in \{ 3,\ldots ,m\}}x_{1}L(E_{1v}^{r_i},\vec{x})=0.
\end{eqnarray*}
And we also change (\ref{eqcit2}) to
\begin{eqnarray*}
\alpha_2(1&-&x_1-x_{t+1}-\cdots -x_v)+\sum\limits_{i\in \{3,\ldots, m\}}L(E_{1\backslash v}^{r_i},\vec{x})+\sum\limits_{v\in \{3,\ldots ,m\}}L(E_{1v}^{r_i},\vec{x})\\
&-&\sum\limits_{i\in \{3,\ldots ,m\}}x_{1}L(E_{1v}^{r_i},\vec{x})=0.
\end{eqnarray*}
We can make other responding changes easily. We omit the detail of the proof here.

By using Theorems \ref{theoPTZ}, \ref{theoPZ}, and \ref{theoTPZZ}, we have

\begin{theorem} \label{th2}
(a) Let  integers $m$ and $t$ satisfy ${t \choose r} \le m \le {t \choose r} + {t-1 \choose r-1}-(2^{r-3}-1)({t-1 \choose r-2}-1)$. Let $H$ be a $\{2,r\}-graph$ ($\{1,2,r\}-graph$, respectively) with $m$ edges in $H^r$ and $t+1$ vertices. If both the vertex order of its maximum complete $\{2,r\}$-subgraphs ($\{1,2,r\}-subgraphs$, respectively) and the vertex order of its maximum complete $\{2\}$-subgraphs ($\{1,2\}-subgraphs$, respectively) are $t$. Then $\lambda'(H) = \lambda'(K_{t}^{\{2,r\}})$ ($\lambda'(H) = \lambda'(K_{t}^{\{1,2,r\}})$, respectively).

(b)  Let  integers $m$ and $t$ satisfy ${t \choose 3} \le m \le {t \choose 3} + {t-1 \choose 2}$. Let $H$ be a $\{1,3\}$-graph ($\{1,2,3\}$-graph, respectively) with $m$ edges in  $H^3$.
If both the order of its maximum complete $\{2,3\}$-subgraphs ($\{1,2,3\}$-subgraphs, respectively) and the order of its maximum complete $\{2\}$-subgraphs ($\{1,2,\}$-subgraphs, respectively) are $t$. Then $\lambda'(H) = \lambda'(K_{t}^{\{2,3\}})$ ($\lambda'(H) = \lambda'(K_{t}^{\{1,2,3\}})$, respectively).

(c)  Let  integers $m$ and $t$ satisfy ${t \choose 3} \le m \le {t \choose 3} + {t-1 \choose 2}-\frac{t}{2}$. Let $H$ be a $\{1,3\}$ with $m$ edges in  $H^3$. Then, if its maximum complete $3$-graph is $K_t^{(3)}$, we have $\lambda'(H)=\lambda'(K_{t}^{\{1,3\}})$;  otherwise
$\lambda'(H)<\lambda'(K_{t}^{\{1,3\}})$.
\end{theorem}

\noindent{\em Proof}
(a) Let $H$ be a $\{2,r\}$-graph with $m$ edges in $H^r$. Let ${\vec x}=(x_1, x_2, \ldots, x_n)$ be an optimal weighting for $H$ for graph-Lagrangian function. Then, use Theorem \ref{MStheo},
\begin{eqnarray*}
\lambda'(H) &=& \lambda'(H,\vec{x})=2!\sum\limits_{ij\in E^{2}}x_ix_j+r!\sum\limits_{\{i_{1}i_{2}\cdots i_{r}\}\in E^{r}}x_{i_{1}} x_{i_{2}}\cdots x_{i_{r}} \\
&\leq&  (1-{1 \over t})+r!\sum\limits_{\{i_{1}i_{2}\cdots i_{r}\}\in E^{r}}x_{i_{1}} x_{i_{2}}\cdots x_{i_{r}}.
\end{eqnarray*}

By Theorem \ref{theoPTZ}, we have $\sum\limits_{\{i_{1}i_{2}\cdots i_{r}\}\in E^{r}}x_{i_{1}} x_{i_{2}}\cdots x_{i_{r}} \leq \lambda'([t]^{(r)})$. Hence $\lambda'(H)\leq (1-{1 \over t})+\lambda([t]^{(r)})=\lambda' (K_{t}^{\{2,r\}})$. On the other side, let $x_1 = x_2 = \ldots =x_{t}=\frac{1}{t}$. We have $\lambda'(H,\vec{x})=\lambda' (K_{t}^{\{2,r\}})$. Therefore $\lambda'(H)=\lambda (K_{t}^{\{2,r\}})$.

The proof of the other results are similar. Note that we use Theorem \ref{theoPZ} in part (b) and Theorem \ref{theoTPZZ} in part (c). We omit the details. \qed

{\bf Acknowledgments} This research is partially  supported by National Natural Science Foundation of China (No. 11271116).

\bibliographystyle{spmpsci}      
\bibliographystyle{unsrt}


\end{document}